\documentclass{amsart}
\usepackage{amsfonts,amsmath,amssymb}

\newtheorem{theorem}{Theorem}[section]
\newtheorem{lemma}[theorem]{Lemma}

\theoremstyle{definition}

\newtheorem{definition}[theorem]{Definition}

\theoremstyle{remark}
\newtheorem{remark}[theorem]{Remark}

\newcommand{\Real}{\mathbb R}
\newcommand{\abs}[1]{\left\vert#1\right\vert}
\newcommand{\set}[1]{\left\{#1\right\}}
\newcommand{\bR}{\mathbb R}
\newcommand{\bS}{\mathbb S}

\DeclareMathOperator{\dist}{dist\,}
\DeclareMathOperator{\diam}{diam} \DeclareMathOperator{\Tr}{tr}
\DeclareMathOperator{\osc}{osc}

\author{Hongjie Dong}
\address{Department of Mathematics, University of Chicago,
5734 S. University Avenue, Chicago, Illinois 60637, USA}
\email{hjdong@math.uchicago.edu}

\author{Seick Kim}
\address{Mathematics Department, University of Missouri,
Columbia, Missouri 65211, USA} \email{seick@math.missouri.edu}

\author{Mikhail Safonov}
\address{School of Mathematics, University of Minnesota,
Minneapolis, Minnesota 55455, USA} \email{safonov@math.umn.edu}

\begin{document}

\title[On Uniqueness of Boundary Blow-up Problem]
{On Uniqueness of Boundary Blow-up Solutions of a Class of
Nonlinear Elliptic Equations}


\begin{abstract}
We study boundary blow-up solutions of semilinear elliptic equations
$Lu=u_+^p$ with $p>1$, or $Lu=e^{au}$ with $a>0$, where $L$ is a
second order elliptic operator with measurable coefficients. Several
uniqueness theorems and an existence theorem are obtained.
\end{abstract}

\maketitle

\section{Introduction}\label{sec:I}

Let $\Omega$ be a bounded domain in $\Real^n$, where $n\ge 2$, and
let $\partial\Omega$ denote its boundary. We consider operators $L$
of the form
\begin{equation*}
L=a^{ij}D_{ij}+b^iD_i-c= a^{ij}(x) \frac{\partial^2}{\partial
x_i\partial x_j}+ b^i(x)\frac{\partial}{\partial x_i}-c(x)
\end{equation*}
whose coefficients $a^{ij} ,b^i, c$ are assumed to be measurable
functions on $\Real^n$ and satisfy
\begin{equation}
                                        \label{eq10.15}
a^{ij}=a^{ji},\quad \sum_i(b^i)^2\leq K,\quad 0\leq c\leq K
\end{equation}
for some fixed constant $K>0$. We also assume that the principal
coefficients $a^{ij}$ satisfy the uniform ellipticity condition;
i.e., there are constants $0<\lambda\le \Lambda<+\infty$ such that
for all $x\in\Omega$, we have
\begin{equation}
\label{eqn:02} \lambda\abs{\xi}^2\le
a^{ij}(x)\xi_i\xi_j\le \Lambda\abs{\xi}^2, \quad  \forall
\xi\in\Real^n.
\end{equation}
Here and throughout the article, the summation convention over
repeated indices is enforced.

In this article, we study the problem
\begin{eqnarray}
\label{eqn:I-2}
&& L u(x)=f(u(x))\quad\text{for}\quad x\in\Omega,\\
\label{eqn:I-3} &&u(x)\to+\infty\quad\mbox{as}\quad d(x):=\dist(x,
\partial\Omega)\to 0,
\end{eqnarray}
where $f(t)=t_{+}^p:=\{\max(t,0)\}^p$ with $p>1$, or $f(t)=e^{at}$.
Solutions of the problem \eqref{eqn:I-2}, \eqref{eqn:I-3} are called
boundary blow-up solutions, or large solutions.

Problems of this type have been studied by many authors.
Bieberbach (1916) considered the equation $\Delta u= e^u$ when
$n=2$, in connection to a problem in Riemannian geometry. Later,
Loewner and Nirenberg (1974) studied the equation $\Delta u
=u_{+}^{(n+2)/(n-2)}$ ($n>2$), which arises in conformal
differential geometry. The problem \eqref{eqn:I-2},
\eqref{eqn:I-3} is also related to probability theory. The
equation $L u=u_{+}^p$, $1<p\le 2$, appears in the analytical
theory of a Markov processes called superdiffusions; see e.g.
Dynkin (2002). By using the potential theory, Labutin (2003)
recently gave a necessary and sufficient Wiener type condition for
the existence of boundary blow-up solutions to $\Delta u=u_+^p$
with $p>1$.

If the domain $\Omega$ is regular enough (e.g. $\Omega$ satisfies
an exterior cone condition), and if the coefficients $a^{ij}$ are
H\"older continuous in $\Omega$, then existence of classical
solutions of the problem \eqref{eqn:I-2}, \eqref{eqn:I-3} can be
established by the method of supersolutions and subsolutions
together with uniform upper bound estimates of Keller (1957) and
Osserman (1957). In fact, Keller and Osserman proved existence of
boundary blow-up solutions of $\Delta u=f(u)$ for a much larger
class of functions $f$ including $f(t)=e^t$ and $f(t)=t^p_{+}$
with $p>1$; see, e.g., the above references for the details.

The question of uniqueness of boundary blow-up solutions has been
studied by many authors. In the case when the domain $\Omega$ is
smooth (e.g. $\Omega$ is of $C^2$), Bandle and Marcus (1992, 1995)
and Lazer and McKenna (1994) proved uniqueness of solutions of the
problem of \eqref{eqn:I-2}, \eqref{eqn:I-3} for a class of
functions $f$ including $f(t)=t^p_{+}$ with $p>1$ by analyzing the
asymptotic behavior of boundary blow-up solutions near the
boundary.

Uniqueness of boundary blow-up solution in non-smooth domains was
also studied by several other authors. Le Gall (1994) investigated
the uniqueness of boundary blow-up solution of $\Delta u =u^2$ in
non-smooth domains by means of a probabilistic representation.
Marcus and V\'eron (1993) proved the uniqueness of boundary
blow-up solution of $\Delta u=u^p$ in very general domains for all
$p>1$, using purely analytical methods. Quite recently, Marcus and
V\'eron (2006) also proved the uniqueness of blow-up solutions for
equation $\Delta u=f(u)$ in bounded domains $\Omega$ such that
$\partial\Omega$ is a locally continuous graph, with convex $f$
satisfying the standard Keller-Osserman condition.

In the main body of this article, we do not impose any regularity
assumptions on the coefficients of operators $L$. In
Theorem~\ref{thm:U-2}, we prove that if $\Omega$ satisfies ``the
uniform exterior ball condition'' (see below for its definition),
then the problem \eqref{eqn:I-2}, \eqref{eqn:I-3} with
$f(t)=t_+^p$ has at most one classical (or strong) solution. A
similar result holds true for $f(t)=e^{at}$ in a special case when
$\Omega$ is convex. Also, in Theorem~\ref{thm:X-2} we show that if
$f(t)=t^p_{+}$ with $p\in(1,1+\frac{2}{\mu(n-1)-1})$, where
$\mu=\Lambda/\lambda\ge 1$, and if $\Omega$ satisfies
$\partial\Omega=\partial\overline\Omega$, then the problem
\eqref{eqn:I-2}, \eqref{eqn:I-3} has at most one classical (or
strong) solution. For the same $f(t)$, by assuming certain
regularity of $a^{ij}$, in Theorem \ref{thm4.2} we prove an
existence and uniqueness result with $p\in (1,1+\frac{2}{n-2})$.
In a special case $L=\Delta$, the results of
Theorems~\ref{thm:X-2} and \ref{thm4.2} are contained in  Marcus
and V\'eron (1997), V\'eron (2001).

Our uniqueness results are based on the iteration technique, which
appears in the proof of Theorem~\ref{thm:U-1}. For operators $L$
with ``good enough'' (e.g. continuous) coefficients, one can also
use another iteration method, introduced by Marcus and V\'eron
(1998), with further development in Marcus and V\'eron (2004). In
particular, they proved (Theorem 3.2 in Marcus and V\'eron, 2004),
that there exists one and only one solution of the problem
\eqref{eqn:I-2}, \eqref{eqn:I-3} in the case
$n=2,\,L=\Delta,\,f(t)=e^t,\,\partial\Omega=\partial\overline\Omega$.
We could not get this result by our method. Roughly speaking, we
need the estimate $e^{u_1}\le Ne^{u_2}$ near $\partial\Omega$ for
any blow-up solutions $u_1$ and $u_2$, while the method in Marcus
and V\'eron (2004) uses a weaker estimate $u_1\le Nu_2$ near
$\partial\Omega$.

The remaining sections are organized in the following way. In
Section~\ref{sec:P}, we give definitions and state some preliminary
lemmas. We state the main results in Section~\ref{sec:M} and prove
them in Sections \ref{proofofthm1} and \ref{sec4}.

\section{Preliminaries}
\label{sec:P}

\begin{definition}
\label{defn:EB} We say that $\Omega$ satisfies the uniform exterior
ball condition with constants $\delta_1\in (0,1)$ and $r_1>0$, if
for arbitrary $x\in\partial\Omega$ and $0<r<r_1$, there exists a
ball $B_{\rho}(y)\subset B_r(x)\setminus\overline\Omega$ with
$\rho=\delta_1 r$.
\end{definition}

\begin{definition}
We say that $u\in C^2(\Omega)$ if $u$ is twice continuously
differentiable in $\Omega$. We write $u\in W^{2,p}_{loc}(\Omega)$
($p\ge 1$) if $u$ is twice weakly differentiable and
$\sum_{\abs{\beta}\le 2} \int_{\Omega'} \abs{D^\beta u}^p <+\infty$
for all $\Omega'\Subset\Omega$. Here $\Omega'\Subset\Omega$ means
that $\Omega'$ is a bounded open set such that its closure
$\overline{\Omega'}$ is a subset of $\Omega$.
\end{definition}

\begin{definition}
We say that $u\in W^{2,n}_{loc}(\Omega)$ is a solution of $Lu=g$ if
$Lu=g$ a.e. in $\Omega$. Similarly, if $u\in W^{2,n}_{loc}(\Omega)$,
then $Lu \ge g$ ($Lu \le g$) in $\Omega$ means $Lu \ge g$ ($Lu \le
g$) a.e. in $\Omega$.
\end{definition}

By Sobolev imbedding theorem, we can always assume that functions in
$W^{2,n}_{loc}(\Omega)$ are continuous in $\Omega$.

\begin{lemma}
        \label{lem:P-1}
Let $\,\Omega\subset\Real^n$ be a bounded domain and let $f$ be an
increasing function. Assume that $u,v\in C^2(\Omega)$ (or $u,v\in
W^{2,n}_{loc}(\Omega)$) satisfy $L u\ge f(u)$ and $L v\le f(v)$ in
$\Omega$. If $\,\liminf_{x\to \partial\Omega}(v-u)(x)\ge 0$, then
$v\ge u$ in $\Omega$.
\end{lemma}

\begin{proof}
Suppose, to the contrary, that there exists $x_0\in\Omega$ such
that $u(x_0)>v(x_0)$. Then for sufficiently small $\epsilon>0$,
$\Omega_\epsilon:=\set{u-v>\epsilon}\ne\emptyset$ and
$\overline\Omega_\epsilon\subset\Omega$. The function
$w:=u-v-\epsilon>0$ in $\Omega_\epsilon$, and $w=0$ on
$\partial\Omega_\epsilon$. Since $f$ is increasing, $L w\ge
f(u)-f(v)\ge 0 $ in $\Omega_\epsilon$. Then, the classical maximum
principle (or Aleksandrov maximum principle) implies $w\le 0$ in
$\Omega_\epsilon$; see e.g. Gilbarg and Trudinger (1983). This
contradiction proves the lemma.
\end{proof}

\begin{lemma}
        \label{lem:P-2}
Let $\Omega\subset\Real^n$ be a bounded domain. If $w\in
C^2(\Omega)$ (or $w\in W^{2,n}_{loc}(\Omega)$) satisfies $L w\ge 0$
in a non-empty subset $\Omega'=\set{x\in\Omega:w(x)>0}$, then
$\partial\Omega'\cap\partial\Omega\neq\emptyset$.
\end{lemma}

\begin{proof}
Otherwise, $\partial\Omega'\subset\Omega$, so that
$\Omega'\Subset\Omega$ and $w=0$ on $\partial\Omega'$. Proceeding
as in the proof of Lemma~\ref{lem:P-1}, we get a contradiction.
\end{proof}

\begin{lemma}
\label{lem:U-5} Let $\Omega=\set{x\in\Real^n:a<\abs{x}<b}$, where
$0\le a<b\le\infty$. If $u(x)=\varphi(\abs{x})$ for some $C^2$
function $\varphi:(a,b)\to \Real$, then, the Hessian $D^2 u(x_0)$
has eigenvalues $\varphi''(r)$ with multiplicity $1$ and
$\varphi'(r)/r$ with multiplicity $n-1$, where $r=\abs{x_0}$.
Therefore, if $\varphi''\ge 0$ and $\varphi'\le 0$, then
\begin{equation}
a^{ij}D_{ij}u(x_0)=\Tr (A\cdot D^2u(x_0)) \ge \lambda \varphi''(r)+
\frac{(n-1)\Lambda}{r} \varphi'(r)
\end{equation}
for any symmetric matrix $A=\{a^{ij}\}$ whose eigenvalues belong to
$[\lambda, \Lambda]$.
\end{lemma}

\begin{proof}
It is a straightforward computation.
\end{proof}

\section{Main results}
                                \label{sec:M}
Our first two results are about the uniqueness of solutions under
the general assumption that the coefficients $a^{ij},b^i,c$ are
measurable functions satisfying \eqref{eq10.15}, \eqref{eqn:02}.

\begin{theorem}
                                \label{thm:U-2}
Let $f(t)=t_{+}^p$ with $p>1$. Assume that $\Omega$ is a bounded
domain satisfying the uniform exterior ball condition. Then there
exists at most one $C^2(\Omega)$ (or $W^{2,n}_{loc}(\Omega)$)
solution of the problem \eqref{eqn:I-2}, \eqref{eqn:I-3}.
\end{theorem}

\begin{theorem}
                                \label{thm3.2}
Let $f(t)=e^{at}$ with $a>0$. Assume that $\Omega$ is a bounded
convex domain in $\Real^n,\,n\ge 2$. Then there exists at most one
$C^2(\Omega)$ (or $W^{2,2}_{loc}(\Omega)$) solution of the problem
\eqref{eqn:I-2}, \eqref{eqn:I-3}.
\end{theorem}

\begin{remark}
                                \label{remark1.16}
Notice that if $u$ is a solution to $Lu=e^{au}$ in $\Omega$, then
$v(x):=au(x/\sqrt{a})$ is a solution to
$$
\bar a^{ij}D_{ij}v+\bar b^iD_{i}v-\bar cv=e^v\quad
\text{in}\,\,\sqrt{a}\,\Omega,
$$
where
$$\bar a^{ij}(x)=a^{ij}(x/\sqrt{a}),\quad \bar
b^i(x)=b^i(x/\sqrt{a})/\sqrt{a},\quad \bar c(x)=c(x/\sqrt{a})/a.
$$ Therefore, without loss of generality, we shall always assume $a=1$ in
the sequel.

\end{remark}

In the next two results, we treat the problem \eqref{eqn:I-2},
\eqref{eqn:I-3} with $f(t)=t^p$ in more general bounded domains
$\Omega\subset \Real^n$. In the case $L=\Delta$, these results are
known from Marcus and V\'eron (1997), V\'eron (2001).

\begin{theorem}
                \label{thm:X-2}
If $p\in(1,1+\frac{2}{\mu(n-1)-1})$, where $\mu=\Lambda/\lambda\ge
1$, and $\partial\Omega=\partial\overline\Omega$, then there exists
at most one $C^2(\Omega)$ (or $W^{2,n}_{loc}(\Omega)$) solution of
the problem \eqref{eqn:I-2}, \eqref{eqn:I-3}.
\end{theorem}

\begin{theorem}
                                            \label{thm4.2}
Suppose that $p\in (1, \infty)$ when $n=2$, and $p\in
(1,\frac{n}{n-2})$ when $n\geq 3$.

i) If $a^{ij}(x)$ are uniformly continuous in a neighborhood of
$\partial \Omega$, and $\partial\Omega=\partial\overline\Omega$,
then there exists at most one solution of the problem
\eqref{eqn:I-2}, \eqref{eqn:I-3}.

ii) If $a^{ij}$ are H\"older continuous in $\Omega$, i.e. $a^{ij}\in
C^\beta(\Omega)$ for some $\beta\in (0,1)$, then there exists at
least one solution of the problem \eqref{eqn:I-2}, \eqref{eqn:I-3}.
\end{theorem}

\begin{remark}
                                    \label{remark3.6}
One can see from the proofs in the following two sections, that the
boundedness assumption of $b^i(x)$ and $c(x)$ can be replaced by
\begin{equation*}
|b^i(x)|=o(d^{-1}(x)), \quad 0\leq c(x)=o(d^{-2}(x)).
\end{equation*}
Also, the uniform ellipticity of the principal coefficients $a^{ij}$
is required only near the boundary $\partial \Omega$, as long as the
weak maximum principle is valid in the entire domain $\Omega$.
Furthermore, if the boundary is smooth (say $C^2$), it suffices to
have $L$ to be nondegenerate only in the normal direction near the
boundary, i.e. there is a $\delta>0$ such that for any $x_0\in
\partial \Omega$ we have $a^{ij}\nu_i\nu_j\geq\lambda$ in
$B_\delta(x_0)$, where $\nu$ is the unit normal direction of
$\partial \Omega$ at $x_0$.
\end{remark}

\begin{remark}
                                        \label{remark3.7}
Without much more work, Theorem \ref{thm:U-2}, \ref{thm3.2} and
\ref{thm:X-2} can be extended to fully nonlinear elliptic equations
$F[u]=u_+^p$ (or $e^{u}$), where $F[u]=F(x,u,Du,D^2u)$ and
$F(x,u,p,q)$ is a function defined on the set
\begin{equation*}
\Gamma:=\bR^n\times \bR \times \bR^n\times \bS^n.
\end{equation*}
Here $\bS^n$ is the set of all symmetric $n\times n$ matrices, and
$F$ satisfies the following natural assumptions
\begin{gather*}
\lambda |\xi |^{2}\leq F(x,u,p,q+\xi \xi ^{T})-F(x,u,p,q)\leq
\Lambda |\xi
|^{2}, \\
|F(x,u,p,q)-F(x,u,p_{1},q)|\leq K|p-p_{1}|, \\
-Ks\leq F(x,u+s,p,q)-F(x,u,p,q)\leq 0,\quad F(x,0,0,0)=0,
\end{gather*}
for any $(x,u,p,q)\in \Gamma$, $s\geq 0$, $p_1\in \bR^n$ and $\xi\in
\bR^n$. In particular, elliptic Bellman equations
$\sup_\beta\{L^\beta u\}=u_+^p$ (or $e^u$) belong to this class,
where linear operators
$$ L^{\beta}=a^{ij}_\beta(x)D_{ij}+b^i_\beta(x)
D_i-c_\beta(x)
$$
satisfy \eqref{eq10.15} and \eqref{eqn:02} with same constants $K,
\lambda, \Lambda$ for all $\beta$. Indeed, it suffices to notice
that under the assumptions above for any two given  $C^2$ functions
$u,v$, we have
$$
F[u]-F[v]=L^{u,v}(u-v),
$$
for some linear  operator $L^{u,v}=a^{ij}D_{ij}+b^iD_{i}-c$
satisfying assumptions \eqref{eq10.15} and \eqref{eqn:02} (see,
for example, Lemma 1.1 in Safonov, 1988). In particular, by
choosing $v\equiv 0$, we get $F[u]=L^u u$ for some linear operator
$L^u$.
\end{remark}

\section{Proof of Theorems \ref{thm:U-2} and \ref{thm3.2}}
                \label{proofofthm1}
Recall the notation $d(x):=\dist(x, \partial \Omega)$. The following
theorem is the main tool of this article in obtaining the uniqueness
results.
\begin{theorem}
\label{thm:U-1} Assume $f(t)=t_{+}^p$ with $p>1$, or $f(t)=e^t$. Let
$\beta=2p/(p-1)$ if $f(t)=t^p_{+}$ with $p>1$, and  $\beta=2$ if
$f(t)=e^t$. If $u_1, u_2$ are $C^2(\Omega)$ (or
$W^{2,n}_{loc}(\Omega)$) solutions of the problem \eqref{eqn:I-2},
\eqref{eqn:I-3}, both satisfying ($i=1, 2$)
\begin{equation} \label{eqn:U-03} N_1 d^{-\beta}\le
f(u_i) \le N_2 d^{-\beta}\quad \text{in}\quad \Delta_\rho:=\{x\in
\Omega\,:\,d(x)< \rho\}
\end{equation}
for some constants $N_1, N_2, \rho>0$, then $u_1\equiv u_2$ in
$\Omega$.
\end{theorem}

\begin{proof}
We first consider the case $f(t)=t^p_{+}$ with $p>1$. Set
$\gamma=2/(p-1)$ so that $\beta=\gamma p=\gamma +2$. Let $u_1$,
$u_2$ be two different $C^2(\Omega)$ (or $W^{2,n}_{loc}(\Omega)$)
solutions of the problem \eqref{eqn:I-2}, \eqref{eqn:I-3}. By
Lemma~\ref{lem:P-1}, they must be different in $\Delta_{\rho}$, and
we may assume that
\begin{equation} \label{u1u2}
u_2(x_0)/u_1(x_0)>k \quad\text{for some}\quad x_0\in \Delta_{\rho}
\quad\text{and}\quad k\ge k_0 >1.
\end{equation}
Note that
$$
L(u_2-ku_1)=f(u_2)-kf(u_1)\geq f(u_2)-f(ku_1)>0\quad\text{in}\quad
\Omega':=\{u_2>ku_1\}.
$$
By Lemma~\ref{lem:P-2} applied to $w=:u_2-ku_1$, we have
$\partial\Omega'\cap\partial\Omega\neq\emptyset$. Therefore, $x_0$
can be chosen arbitrary close to $\partial\Omega$, and we may assume
that
\begin{equation}
                                        \label{eq11.38}
B_{r}(x_0)\subset \Delta_\rho\quad\text{and}\quad K(2r+r^2)\leq
\Lambda, \quad\text{where } r:=d(x_0)/2.
\end{equation}
The set $\Omega_0:=\set{u_2-ku_1>0}\cap B_r(x_0)\Subset\Omega$. In
$\Omega_0$, we have $r<d(x)<3r$, and
\begin{equation*}
L(u_2-k u_1)=u_2^p-ku_1^p>(k^p-k)u_1^p\ge (k_0^{p-1}-1)ku_1^p \ge
c_0 kr^{-\beta},
\end{equation*}
where $c_0:=(k_0^{p-1}-1)3^{-\beta}N_1>0$. On the other hand, the
function
\begin{equation}
\label{eqn:U-02} w(x)=c_1 k r^{-\beta}(r^2-\abs{x-x_0}^2),
\quad\text{where } c_1:=c_0/(3n\Lambda),
\end{equation}
satisfies
$$ Lw\ge -c_1 k r^{-\beta} (2n\Lambda +2Kr+Kr^2) \ge -c_0 kr^{-\beta}
\quad\text{in}\quad B_r(x_0)\supset \Omega_0.
$$
Then the function $w_1:=u_2-ku_1+w$ satisfies $Lw\ge 0$ in
$\Omega_0$, and by the maximum principle, it attains its maximum on
$\overline{\Omega_0}$ at some point $x_1\in \partial\Omega_0$. Note
that $x_1$ cannot belong to $B_r(x_0)$, because on the set
$(\partial\Omega_0)\cap B_r(x_0)$, we must have $u_2=ku_1$, which in
turn implies $w_1=w\le w(x_0) <w_1(x_0)\le w_1(x_1)$. Therefore,
$x_1\in \partial B_r(x_0)$, so that $w(x_1)=0$, and
\begin{equation} \label{eqn:U-04}
(u_2-k u_1)(x_1) = w_1(x_1)\ge w_1(x_0) >w(x_0) = c_1 k
r^{2-\beta}=c_1kr^{-\gamma}.
\end{equation}

Since $d(x_1)\ge r$, from \eqref{eqn:U-03} it follows
\begin{equation*}
    f(u_1(x_1))=u_1^p(x_1)\le N_2 d^{-\beta}(x_1) \le N_2 r^{-\beta}=N_2 r^{-\gamma p}.
\end{equation*}
This estimate together with \eqref{eqn:U-02}, \eqref{eqn:U-04}
implies
\begin{equation*}
u_2(x_1)>(1+c_2) k u(x_1),\quad \text{where } c_2:= c_1
N_2^{-1/p}>0.
\end{equation*}
Again, replacing $x_1$ by another point near $\partial\Omega$ if
necessary, we may assume that \eqref{u1u2}, \eqref{eq11.38} hold
with $(1+c_1)k,x_1,r_1:=d(x_1)/2$ in place of $k,x_0,r$
respectively. By iterating, we obtain a sequence
$\set{x_j}_{j=0}^\infty\subset \Delta_\rho$ such that
$u_2(x_j)/u_1(x_j)>(1+c_2)^j k_0$, which tends to infinity. However,
\eqref{eqn:U-03} implies $u_2/u_1\le (N_2/N_1)^{1/p}$ in
$\Delta_\rho$. This contradiction proves that $u_1\equiv u_2$ in
$\Omega$.

Now we consider the case $f(t)=e^t$. We proceed similarly as above.
Let $u_1,u_2$ be two different solutions of the problem
\eqref{eqn:I-2}, \eqref{eqn:I-3}. We may assume that
\begin{equation*}
u_2(x_0)-u_1(x_0)>k\quad\text{for some}\quad x_0\in
\Delta_{\rho}\quad\text{and}\quad k\ge k_0>0.
\end{equation*}
By Lemma~\ref{lem:P-2} applied to $w:=u_2-u_1-k$, we may also assume
that $x_0$ is chosen such that \eqref{eq11.38} holds. Then
\begin{equation*}
L(u_2-u_1-k)\ge e^{u_2}-e^{u_1}>(e^k-1)e^{u_1}\ge c_3k r^{-2}
\end{equation*}
on the set $\Omega_0:=\{u_2-u_1>k\}\cap B_r(x_0)$, where
$c_3:=N_1/9$. On the other hand, the function
\begin{equation*}
w(x):=c_4 k r^{-2}(r^2-\abs{x-x_0}^2), \quad\text{where} \quad
 c_4:=c_3/(3n\Lambda)>0 ,
\end{equation*}
satisfies $Lw\ge -c_3 kr^{-2}$ in $\Omega_0$. Then the function
$w_1:=u_2-u_1-k+w$ satisfies $Lw_1\ge 0$ in $\Omega_0$, hence it
attains its maximum on $\overline{\Omega_0}$ at some point
$x_1\in\partial\Omega_0$, which cannot belong to $B_r(x_0)$.
Therefore, $x_1\in \partial B_r(x_0),\,w(x_1)=0$, and
\begin{equation*}
    u_2(x_1)-u_1(x_1)-k=w_1(x_1)\ge w_1(x_0)> w(x_0)=c_4 k.
\end{equation*}
As before, by iterating this process, we obtain a sequence
$\set{x_j}_{j=0}^\infty\subset \Delta_\rho$ such that
$u_2(x_j)-u_1(x_j)>(1+c_4)^j k_0$, which tends to infinity.
However, \eqref{eqn:U-03} implies that $u_2-u_1\le \ln(N_2/N_1)$
in $\Delta_\rho$. Again, this contradiction leads to the
conclusion that $u_1\equiv u_2$ in $\Omega$. The theorem is
proved.
\end{proof}

\begin{remark}
By easy modifications of the proof above, one can see that Theorem
\ref{thm:U-1} can be extended to any locally Lipschitz, increasing
function $f$, which is equal to $e^t$ in $(N_1, \infty)$, or $f$
which is equal to $t^p$ in $(N_1, \infty)$, for some $N_1>0$, and
satisfies the additional condition $f(\mu t)\geq \mu f(t)$ for any
$\mu\ge 1$ and $t\in \bR$.
\end{remark}

We derive a lower and upper bounds in the following two lemmas.
\begin{lemma}
\label{lem:U-1} Let $f(t)$, $\beta$, and $\Delta_\rho$ be as in
Theorem~\ref{thm:U-1}. If $u$ is a $C^{2}(\Omega)$ (or
$W^{2,n}_{loc}(\Omega)$) solution of the problem \eqref{eqn:I-2},
\eqref{eqn:I-3}, then we have
\begin{equation}
                                \label{eqn:U-13}
1/N_1 \le f(u) \,\,\text{in}\,\,\Omega,\quad f(u) \le N_2 d^{-\beta}
\,\, \text{in}\,\, \Delta_1:=\{x\in\Omega:\,d(x)<1\}.
\end{equation}
Here $N_1, N_2$ are positive constants depending only on $\,n,
\lambda, \Lambda, K,$ and $p$ if $f(t)=t_{+}^p$ with $p>1$; $N_1$
may also depend on $\diam\Omega$.
\end{lemma}

\begin{proof}
First, we consider the case $f(t)=t_{+}^p$ with $p>1$. Without loss
of generality, we assume that $\Omega$ contains the origin. The
lower bound follows from an observation that $\varepsilon e^{\eta
x_1}$ is a bounded subsolution if we first choose $\eta$
sufficiently large, and then $\varepsilon>0$ sufficiently small. It
remains to get the upper bound.  Fix $x_0\in\Delta_1$ and
$r<d(x_0)<1$. Denote
\begin{equation*}
w_0(x):= N_0(1-\abs{x}^2)^{-\gamma},\quad
w(x):=r^{-\gamma}w_0((x-x_0)/r),
\end{equation*}
where $\gamma:=2/(p-1)$ as before. If we set $N_0:=(2\gamma
(n+2\gamma)\Lambda+2 \gamma K)^{\gamma/2}$, then we have $L w \le
w^p$ in $B_r(x_0)$ for any elliptic operator $L$ whose coefficients
satisfy \eqref{eq10.15}, \eqref{eqn:02}. By Lemma~\ref{lem:P-1},
$u(x)\le w(x)$ in $B_r(x_0)$. In particular, we have $u(x_0)\le
w(x_0)=N_0r^{-\gamma}$. Therefore, we get the desired bound
\eqref{eqn:U-13} with $N_2:=N_0^p$ by letting $r\to d(x_0)$.

The case $f(t)=e^t$ is treated similarly. Without loss of
generality, we may assume that $\Omega$ lies in the half-space
$\{x_1>0\}$. Fix positive constants $\eta_1$ and $\eta_2$, such that
$$\lambda\eta_1^2-K\eta_1-K\ge 1,
\quad \eta_2\ge\sup_{\Omega}e^{\eta_1 x_1}.$$ Then the function
$v:=e^{\eta_1x_1}-\eta_2$  satisfies $v\le 0$ and
$$Lv=(a^{11}\eta_1^2 +b^1\eta_1 -c)e^{\eta_1 x_1}+c\eta_2
\ge (\lambda\eta_1^2-K\eta_1-K)e^{\eta_1 x_1}\ge 1\ge e^v$$ in
$\Omega$. Hence $u\ge v$ in $\Omega$, and the lower bound follows.
For the upper bound, we fix $x_0\in\Delta_1$ and set
\begin{equation*}
w_0(x):=\ln N_2-2 \ln(1-\abs{x}^2),\quad w(x):=w_0((x-x_0)/r)-2\ln
r,
\end{equation*}
where $N_2:= 4n(\Lambda+K)$. Then $Lw\le e^w$ in $B_r(x_0)$.
Again, Lemma~\ref{lem:P-1} implies that $u(x_0)\le w(x_0)=
\ln(N_2/r^2)$. By letting $r\to d(x_0)$, we obtain the bound
\eqref{eqn:U-13}. The lemma is proved.
\end{proof}

\begin{remark}
\label{rmk:3} In the previous lemma, the assumption \eqref{eqn:I-3}
was used only for the proof of the lower bound in \eqref{eqn:U-13}.
Note that the upper bound
\begin{equation}
\label{eqn:13} f(u(x)) \le N_2 d^{-\beta}(x)\qquad \forall
x\in\Omega
\end{equation}
is valid for any $C^2(\Omega)$ (or $W^{2,n}_{loc}(\Omega)$) solution
$u$ of \eqref{eqn:I-2}.
\end{remark}

\begin{lemma}
\label{lem:U-2} Let $\Omega$ be a bounded domain satisfying the
uniform exterior ball condition with constants $r_1$ and $\delta_1$
(see Definition~\ref{defn:EB}). Assume $f(t)=t_+^p$, where $p>1$,
and set $\beta:=2p/(p-1)$. If $u$ is a $C^2(\Omega)$ (or
$W^{2,n}_{loc}(\Omega)$) solution of the problem \eqref{eqn:I-2},
\eqref{eqn:I-3}, then
\begin{equation}
\label{eqn:15} f(u)\ge N
d^{-\beta}\quad\text{in}\quad\Delta_{\rho}:=\{x\in\Omega:\,d(x)<\rho\},
\end{equation}
where $\rho:=\min(r_1,1/2)$, and  $N>0$ is a constant depending only
on $n, \lambda, \Lambda, K, p, \delta_1$, and  $r_1$.
\end{lemma}
\begin{proof}
For a fixed point $x_0\in\Delta_{\rho}$, choose
$z_0\in\partial\Omega$ such that $|x_0-z_0|=r_0:=d(x_0)$, and then
$y_0$ such that $B_{\delta_1 r_0}(y_0)\subset B_{r_0}(z_0)\setminus
\overline\Omega$.

Set $\delta:=\delta_1/2$ and $r:=2r_0$. Observe that if $m=m(K,n,
\delta)$ is sufficiently large, then $v_0(t):=(1-t)^m$ satisfies
\begin{equation}
\label{eqn:17} \left\{
\begin{aligned}
&\lambda v_0''(t) + \frac{(n-1)\Lambda}{t}\,v_0'(t)+ K v_0'(t)- K
v_0(t) \ge v_0^p(t),
\quad \forall t\in(\delta, 1), \\
& v'_0(t) < 0, \quad \forall t\in(\delta, 1),\\
& v_0(1)=0,\quad v_0(t)>0 \quad \forall t\in[\delta,1).
\end{aligned}
\right.
\end{equation}
Note that $\delta r=\delta_1 r_0<r=2r_0 <2\rho\le 1$. Using
Lemma~\ref{lem:U-5}, it is easy to check that the function
$v(x):=r^{-\gamma}v_0(\abs{x-y_0}/r)$ with $\gamma:=2/(p-1)$
satisfies $Lv \ge v^p$ in $\Omega\cap B_{r}(y_0)$. Thus, by
Lemma~\ref{lem:P-1}, since $|x_0-y_0|\le
(2-\delta_1)r_0=(1-\delta)r$, we have
\begin{equation*}
u(x_0)\ge v(x_0)\ge 2^{-\gamma} d^{-\gamma}(x_0)\, v_0(1-\delta)
\end{equation*}
From here the desired lower bound follows with $N:= 2^{-\gamma
p}v_0^p(1-\delta)$. The lemma is proved.
\end{proof}

Now we are ready to prove Theorems \ref{thm:U-2} and \ref{thm3.2}.

{\bf Proof of Theorem \ref{thm:U-2}:} It follows readily from
Theorem~\ref{thm:U-1}, Lemma~\ref{lem:U-1}, and Lemma~\ref{lem:U-2}.

{\bf Proof of Theorem \ref{thm3.2}:} Fix a constant
$D>\diam(\Omega)$. We may assume that
$$\overline{\Omega}\subset \{x=(x_1, \ldots,x_n)\in
\Real^n:\,\,0<x_1<D\}.$$ Note that the function $v_0:=-2\ln x_1$
satisfies
$$Lv_0=2a_{11}x_1^{-2}-2b^1x_1^{-1}+2c\ln x_1\ge
2\lambda x_1^{-2}-2Kx_1^{-1}-2K |\ln x_1|$$ for $x_1>0$. Choose
constants $\delta=\delta(\lambda,K)\in (0,1)$ and
$N=N(\lambda,K,D)\ge K\,|\ln\lambda|$, such that
\begin{eqnarray*}
  Lv_0 \ge \lambda x_1^{-2} &\quad \text{for}\quad & 0<x_1<\delta, \\
  N-c\ln\lambda+Lv_0 \ge \lambda x_1^{-2} & \quad\text{for}\quad & \delta\le
  x_1<D.
\end{eqnarray*}
As in the proof of Lemma~\ref{lem:U-1}, take a function
$v:=e^{\eta_1x_1}-\eta_2$ satisfying
$$v\le 0,\quad Lv\ge 1\quad\text{for}\quad 0<x_1<D.$$
Then the function $w:=Nv+\ln\lambda +v_0$ satisfies
$$Lw\ge N-c\ln\lambda +Lv_0\ge \lambda x_1^{-2}
=e^{\ln\lambda+v_0} \ge e^w \quad\text{for}\quad 0<x_1<D.$$ By
Lemma~\ref{lem:P-1}, we must have
$$u\ge w,\quad e^u\ge e^w \ge \lambda\,e^{-N\eta_2} x_1^{-2}
\quad\text{in}\quad\Omega.$$

Finally, note that the conditions \eqref{eq10.15}, \eqref{eqn:02} on
the coefficients of $L$ are invariant with respect to parallel
translations and rotations in $\Real^n$. Therefore, for any fixed
$x=(x_1, \ldots,x_n)\in\Omega$, we can always assume that
$x_2=\cdots=x_n=0$, and $x_1>0$ can be made arbitraryly close to
$d(x)=\dist(x, \partial\Omega)$. This means that we have the lower
bound
$$e^u \ge \lambda\,e^{-N\eta_2} d^{-2}=:N_1 d^{-2}
\quad\text{in}\quad\Omega.$$ This estimate, together with Theorem
\ref{thm:U-1} and the upper bound in Lemma \ref{lem:U-1}, yields
the uniqueness.

\section{Proof of Theorems \ref{thm:X-2} and \ref{thm4.2} }
                                            \label{sec4}

In this section, we prove the uniqueness of a solution of the
problem \eqref{eqn:I-2}, \eqref{eqn:I-3} with $f(t)=t_+^p$, in more
general domains.

{\bf Proof of Theorem \ref{thm:X-2}:} Assume that $u$ is a
$C^2(\Omega)$ (or $W^{2,n}_{loc}(\Omega)$) solution of the problem
\eqref{eqn:I-2}, \eqref{eqn:I-3}, and set
\begin{equation*}
    \gamma:=2/(p-1),\quad
    \gamma_0:=\gamma((\gamma+1)\lambda+(1-n)\Lambda).
\end{equation*}
Note that the assumption $p\in(1,1+\frac{2}{\mu(n-1)-1})$ implies
$\gamma_0>0$. Let $r_0\in (0,1)$ be such that $2(Kr_0^2+K\gamma
r_0)\leq \gamma_0$. Fix $x_0\in\Delta_{r_0/2}$ and choose
$z_0\in\partial\Omega$ such that $|x_0-z_0|=r:=d(x_0)$. From
$\partial\Omega=\partial\overline\Omega$ it follows that there
exists a point $y_0\in B_{r/2}(z_0)\setminus\overline\Omega$. Using
Lemma~\ref{lem:U-5}, it is easy to check that the function
\begin{equation*}
    v(x):=c_0\abs{x-y_0}^{-\gamma}-c_0(2r)^{-\gamma},\quad\text{where}\quad
    c_0:=(\gamma_0/2)^{\gamma/2},
\end{equation*}
satisfies $Lv\ge v^p$ in $\Omega\cap B_{2r}(y_0)$. Moreover, $v\in
C^2(\overline\Omega),\,v < +\infty$ on $(\partial\Omega)\cap
\overline B_{2r}(y_0)$, and $v=0$ on $\overline \Omega\cap \partial
B_{2r}(y_0)$. Therefore, by Lemma~\ref{lem:P-1}, $u(x)\ge v(x)$ in
$\Omega$. In particular, we have
$$
u(x_0)\ge v(x_0)\ge c_1 d^{-\gamma}(x_0),\quad\text{where}\quad
c_1:=(1.5^{-\gamma}-2^{-\gamma})c_0.$$ Also, by Lemma~\ref{lem:U-1},
$u(x_0)\le c_2 d^{-\gamma}(x_0)$ in $\Delta_1$, for some $c_2>0$
depending only on $n, \lambda, \Lambda, K$, and $p$. Since $x_0\in
\Delta_{r_0/2}$ is arbitrary, we have proved that
\begin{equation*}
c_1 d^{-\gamma} \le u \le c_2 d^{-\gamma} \quad
\text{in}\,\,\Delta_{r_0/2}.
\end{equation*}
Now the desired statement follows from Theorem~\ref{thm:U-1}.
\medskip

{\bf Proof of Theorem \ref{thm4.2}:} We prove part ii) first. Let
$\set{\Omega_m}_{m=1}^\infty$ be an exhausting sequence of smooth
subdomains of $\Omega$; i.e.,
$\Omega_m\Subset\Omega_{m+1}\Subset\Omega$ and
$\bigcup_{m=1}^\infty\Omega_m=\Omega$. Let $u_m$ be the unique
boundary blow-up solution of $L u = u^p_{+}$ in $\Omega_m$ for
each $m\ge 1$. (For existence of such solutions $u_m$, see, e.g.
Keller, 1957; uniqueness is a consequence of
Theorem~\ref{thm:U-2}.) By Lemma~\ref{lem:P-1},
$\set{u_m}_{m=1}^\infty$ is a decreasing sequence, and by
Lemma~\ref{lem:U-1}, it is bounded below by some constant
$1/N_1>0$. Hence, the limit function $u$ exists in $\Omega$ and by
the standard elliptic theory, it is a solution of $Lu=u^p$ in
$\Omega$.

We claim that $u$ is indeed a boundary blow-up solution. In order to
prove this, it suffices to show that for any $y_0\in \partial
\Omega$,
\begin{equation}
                                            \label{eq5.33}
u_m(x)\geq N_0|x-y_0|^{-\gamma}\quad \text{in}\quad \Omega_m\cap
B_{r_0}(y_0),
\end{equation}
where $\gamma=2/(p-1)$ as before, and $N_0,r_0$ are positive
constants independent of $m$.

We first do a linear transformation to make $y_0=0$,
$a^{ij}(y_0)=\delta^{ij}$, and still use the same notations for
simplicity. Due to \eqref{eqn:02}, the scales in these two
coordinate systems are comparable. Therefore, we only need to verify
\eqref{eq5.33} in the new coordinates. Set
\begin{equation*}
    v_0(x)=c_p\abs{x}^{-\gamma},\quad\text{where}\quad
    c_p:=\{\gamma(\gamma+2-n)/2\}^{\gamma/2}.
\end{equation*}
Since $a^{ij}$ are uniformly continuous, one can choose $r_1>0$
sufficiently small, such that in $\Omega_m\cap B_{r_1}(0)$,
\begin{eqnarray*}
Lv_0(x)&=& c_p\left\{\Delta(|x|^{-\gamma})+(a^{ij}(x)-a^{ij}(0))
D_{ij}(|x|^{-\gamma})\right.\\
&+&\left. b^iD_i (|x|^{-\gamma})-c|x|^{-\gamma}\right\}
\\
&\geq& c_p\left\{\gamma(\gamma+2-n)-K\gamma r_1-Kr_1^2+
N(n,p)\,\omega(r_1)\right\}|x|^{-\gamma-2}
\\
&\geq& (c_p/2)\gamma(\gamma+2-n)|x|^{-\gamma-2} =v_0^p(x),
\end{eqnarray*}
where $\omega(r_1)= \max_{i,j} \{\osc_{\Omega\cap
B_{r_1}(0)}a^{ij}\}$. Then the function
$v(x):=c_p|x|^{-\gamma}-c_pr_1^{-\gamma}$ satisfies $Lv \ge v^p$ in
$\Omega_m\cap B_{r_1}(0)$, and $v(x)=0$ on $\partial B_{r_1}(0)$.
Therefore, by Lemma~\ref{lem:P-1}, we have $u_m(x)\ge v(x)$ in
$\Omega_m\cap B_{r_1}(0)$, and the desired estimate \eqref{eq5.33}
follows with $N_0:=c_p(1-2^{-\gamma})$ and $r_0:=r_1/2$ .

For the proof of  i), due to Theorem \ref{thm:U-1}, it suffices to
get the estimate
$$
N_1d^{-\gamma}\leq u\leq N_2d^{-\gamma}
$$ in a neighborhood of $\partial \Omega$. Arguing as in the proof of
Theorem \ref{thm:X-2}, this estimate can be proved by using the
barrier function $v(x)$ constructed in the proof of ii). The
details are left to the reader.

Acknowledgements The authors are thankful to Laurent V\'eron for
very useful discussion.

\end{document}